\documentclass[11pt,reqno]{amsart}
\usepackage{amssymb,amsthm,amsfonts,amsmath,amscd,amssymb,epsfig}

\usepackage{enumerate,array,hyperref, graphicx,color,amsthm,
bbm,mathrsfs, geometry}

\theoremstyle{plain}
\newtheorem{theorem}{Theorem}[section]

\newtheorem{proposition}[theorem]{Proposition}

\theoremstyle{remark}

\newtheorem*{question}{Question}
\newtheorem*{remark}{Remark}

\theoremstyle{definition}
\newtheorem*{definition}{Definition}
\newcommand{\R}{{\mathbb{R}}}
\newcommand{\C}{{\mathbb{C}}}

\renewcommand{\min}{{\rm min}}

\newcommand{\comment}[1]{}

\newcommand{\beq}{\begin{equation}}
\newcommand{\beqn}{\begin{equation}\nonumber}
\newcommand{\eeq}{\end{equation}}

\newcommand{\bea}{\begin{equation}\begin{aligned}}
\newcommand{\bean}{\begin{equation}\begin{aligned}\nonumber}
\newcommand{\eea}{\end{aligned}\end{equation}}

\renewcommand{\Re}{\,\mathrm{Re}}
\renewcommand{\Im}{\,\mathrm{Im}}

\begin{document}

\title{ Global surfaces of section in the planar restricted 3-body problem}
\author{Peter Albers}
\author{Joel W.~Fish}
\author{Urs Frauenfelder}
\author{Helmut Hofer}
\author{Otto van Koert}
\address{
	Peter Albers\\
	Department of Mathematics\\
	Purdue University
	}
\email{palbers@math.purdue.edu}
\address{
	Joel W.~Fish\\
	Department of Mathematics\\
	Stanford University
	}
\email{joelfish@math.stanford.edu}
\address{
	Helmut Hofer\\
	Institute for Advanced Study
	}
\email{hofer@ias.edu}
\address{
    Urs Frauenfelder\\
    Department of Mathematics and Research Institute of Mathematics\\
    Seoul National University}
\email{frauenf@snu.ac.kr}
\address{
    Otto van Koert\\
    Department of Mathematics and Research Institute of Mathematics\\
    Seoul National University}
\email{okoert@snu.ac.kr}

\keywords{Restricted 3-body problem, global surface of section, convexity}

\maketitle

\begin{abstract}
The restricted planar three-body problem has a rich history, yet many unanswered questions still remain.
In the present paper we prove the existence of a global surface of section near the smaller body in a new  range of energies and mass ratios
for which the Hill's region still has three connected components. The
approach relies on recent global methods in symplectic geometry
and contrasts sharply with the perturbative methods used until now.

\end{abstract}

%\tableofcontents

%

\section{Introduction}

We consider the planar restricted 3-body problem in rotating coordinates. We call the two primaries the sun and the earth and place the earth at the origin of our coordinate system. Then the describing Hamiltonian $H:\C\setminus\{0,1\}\times\C\to\R$ is given by
\beq
H(q,p)=\frac12|p|^2+\langle p,iq\rangle -\langle p,i\mu\rangle-\frac{1-\mu}{|q|}-\frac{\mu}{|q-1|}\;,
\eeq
where $\langle p,iq\rangle=p_2q_1 - p_1 q_2$. Also $\mu\in[0,1]$ is the mass ratio $\mu=\frac{m_S}{m_E+m_S}$ where $m_E$ is the mass of the earth and $m_S$ is the mass of the sun. The approximate  value
for the (real) sun/earth system is $\mu\approx0.999997$. For $\mu\notin \{0,1\}$ the Hamiltonian $H$ has five critical points $L_1,\ldots,L_5$, which we order by increasing value of $H$. These are the Lagrange points.  For $\mu\in \{0,1\}$ the Hamiltonian $H$ has a critical manifold diffeomorphic to $S^1$. Note that the critical value $H(L_1)$ converges to $-\frac32$ as $\mu$ tends to either $0$ or $1$.

If we choose the energy level $-c$ to be below the first Lagrange value $H(L_1)$, then the energy hypersurface $H^{-1}(-c)$ has three connected components: one is near the earth, one is near the sun, and one is near infinity. Throughout this article, we will focus on the component closest to the earth. Note that the components around the earth and sun are non-compact due to collisions with the respective primaries.
However, it is well known that such two-body collisions can be regularized. Recall the Levi-Civita coordinates given by
$q=2v^2$ and $p=\frac{u}{\bar{v}}$ in \cite{Levi_Civita}.  These coordinates define  a 2:1-map, which is symplectic up to a factor $4$.  Indeed, $\Re(dq\wedge d\bar{p})=4\Re(dv\wedge d\bar{u})$. Transforming and regularizing the Hamiltonian function at energy $-c$ leads to
\beq
K_{\mu,c}(v,u):=|v|^2\big(H(v,u)+c\big)=\frac12|u|^2+2|v|^2\langle u,iv\rangle -\mu\Im(uv)-\frac{1-\mu}{2}-\frac{\mu|v|^2}{|2v^2-1|}+c|v|^2\;.
\eeq
For $\mu\notin\{0,1\}$ the component of the energy hypersurface $H^{-1}(c)$ around the earth lifts to a compact component $\Sigma_{\mu,c}$ of the energy hypersurface $K_{\mu,c}^{-1}(0)$.
The energy surface  $K_{\mu,c}^{-1}(0)$ is diffeomorphic to $S^3$.

Next we recall a version of the definition of a surface of section.
\begin{definition}
Let $\Sigma$ be a smooth closed three-manifold equipped with a smooth flow without rest points.
A \emph{global disk-like surface of section} consists of a topologically embedded  closed  disk ${\mathcal D}\subset \Sigma$ having the following properties:
\begin{itemize}
\item[(1)]  The boundary $\partial{\mathcal D}$ is an (un-parametrized) periodic orbit, called the spanning orbit.
\item[(2)] The interior of the disk $\dot{\mathcal D}={\mathcal D}\setminus\partial{\mathcal D}$ is a smooth submanifold
of $\Sigma$ and is transversal to the flow.
\item[(3)] Every orbit, other than the spanning orbit, intersects the (interior of the) disk in forward and backward time.
\end{itemize}
\end{definition}
The above definition allows the disk ${\mathcal D}$ to be rather wild near its boundary. Given a global disk-like surface of section
it follows that there exists a smooth map $\psi:\dot{\mathcal D}\rightarrow\dot{\mathcal D}$, called the global return map.
In general $\psi$, which is defined on the interior of the disk, does not need to have an
extension to the boundary. Note that there is not much one can say about a continuous self-map defined on an open disk. For example, Brouwer's fixed point theorem fails.
However, much more can be said if the map is an area-preserving diffeomorphism. Indeed, a consequence of Brouwer's translation theorem is that such maps always have a fixed point.

The notion of global surfaces of section goes back to Poincar\'e, and it is clear that they encode much of the dynamics on the energy surface. Later we shall describe some consequences of their existence, but presently we state our main result.

\begin{theorem}\label{thm:main}
For every $c>\frac32$ there exists $\mu_0=\mu_0(c)\in [0,1) $ such that for all $\mu_0<\mu<1$ there exists a   global disk-like surface of section for the component $\Sigma_{\mu,c}$ of the energy hypersurface of $K_{\mu,c}^{-1}(0)$.
\end{theorem}

The existence of this global surface of section follows as a consequence of a global result in symplectic geometry, which is applicable provided
the energy surface satisfies certain geometric conditions. Also note that it seems impossible to obtain this surface of section by the usual method of perturbing an understood model. Instead, to achieve our result we must verify that a certain convexity assumption holds.  We make this precise below.

\begin{definition}
The \emph{convexity range} $\mathfrak{C}$ is defined to be the collection of all pairs $(c,\mu)$ with $c>\frac{3}{2}$ and $\mu\in (0,1)$ such that
the energy surface $\Sigma_{\mu,c}$ bounds a strongly convex domain.
\end{definition}

Here we say a compact surface $\Sigma\subset \R^4$ bounds a strongly convex domain provided there exists a constant $\delta>0$ and a smooth convex function $C:{\mathbb R}^4\to\R$ such that $\Sigma=C^{-1}(1)$ and the matrix valued function $D^2C(z)- \delta Id$ is positive definite for all $z\in \R^4$. An elementary exercise shows that a connected compact hypersurface  $\Sigma\subset {\mathbb R}^{4}$ bounds a strongly convex domain $W$ whenever there exists a smooth function $\phi:{\mathbb R}^4\rightarrow {\mathbb R}$ with the following properties:
\begin{enumerate}
\item $\Sigma=\{\phi=0\}$ is a regular level of $\phi$.
\item $W=\{z\in \R^4:\phi(z)\leq 0\}$ is bounded.
\item $D^2\phi(z)(h,h)>0$ for each point $z\in W$ and for each non-zero vector $h$.
\end{enumerate}

We now state the main technical result of this article.

\begin{proposition}\label{prop1}
For each $c>\frac{3}{2}$ there exists a number $\mu_0(c)\in (0,1)$ such that
$$
\{(c,\mu)\ |\ c>{\textstyle \frac{3}{2}},\ \mu\in(\mu_0(c),1)\}\subset \mathfrak{C}.
$$
\end{proposition}
As the  elementary proof of Proposition \ref{prop1} shows,
it should be possible to use a computer to get a more precise idea of the convexity range.

Observe that Theorem \ref{thm:main} now follows from Proposition \ref{prop1} and the following theorem which relies on a pseudoholomorphic curve theory for contact manifolds.
 The core idea is the construction of certain foliations, called finite energy foliations.
\begin{theorem}[\cite{HWZ_the_dynamics_on_three_dimensional_strictly_conve_eneergy_surfaces}]\label{thm2}
If $\Sigma$ is a smooth, regular, bounded energy surface in ${\mathbb R}^4$ bounding a strongly convex domain,
then there exists a global disk-like surface of section ${\mathcal D}$ and an associated global return map
$\psi:\dot{\mathcal D}\rightarrow \dot{\mathcal D}$, which is smoothly conjugated to a smooth area-preserving
disk map $\Psi:\dot{D}\rightarrow \dot{D}$, where $\dot{D}$ is the open unit disk in the plane equipped with the Lebesgue measure.
\end{theorem}

A celebrated result by Franks, \cite{Franks_Geodesics_on_S2_and_periodic_points_of_annulus_homeomorphisms}, implies that $\Psi$ either has precisely one periodic point or infinitely many.
This result then also holds for $\Sigma_{\mu,c}$ whenever $(c,\mu)\in \mathfrak{C}$.

\begin{remark}
Analyzing the proof of Theorem \ref{thm:main} in \cite{HWZ_the_dynamics_on_three_dimensional_strictly_conve_eneergy_surfaces} one should
be able to obtain some refinements. First, one should be able to find a continuously differentiable ${\mathcal D}$ in such a way that the return map
defined on $\dot{\mathcal D}$ has a continuously differentiable extension over the closed disk. This map should be conjugated to an area-preserving map $\Psi$ on the closed unit disk $D$. Recent results by Franks/Handel,\cite{Franks_Handel_Periodic_points_of_Hamiltonian_surface_diffeomorphisms}, and LeCalvez, \cite{LeCalvez_Periodic_orbits_of_Hamiltonian_homeomorphisms_of_surfaces}, then imply that
for $ \phi$ one of the following holds:
 \begin{itemize}
 \item[(1)]  $\Psi$ is a pseudo-rotation, i.e. it has precisely one periodic point.
 \item[(2)] Some iterate of $\Psi$ is the identity.
 \item[(3)] The minimal periods of periodic orbits of $\Psi$  are unbounded.
  \end{itemize}
As shown in \cite{Bramham1,Bramham2,Bramham3} finite energy foliations can also be used to study area-preserving disk maps.
 We refer the reader for more details to \cite{Bramham_Hofer_First_steps_toward_a_symplectic_dynamics} where some of the recent results on area-preserving disk maps are surveyed.
 We leave the construction described in the remark to the interested reader.
 One should be able to prove that item (2) would imply integrability of the flow on the corresponding energy surface $\Sigma_{\mu,c}$, which seems unlikely.
 Also item (1) seems unlikely for most energy surfaces.
Finally, observe that in the simplest case, namely the rotating Kepler problem, for which $\mu=0$, item (3) holds.
\end{remark}
{\bf Acknowledgments:} We thank B.~Bramham and E.~Belbruno for stimulating discussions.
The research of  P.~Albers, J.~Fish and H.~Hofer  was partially supported by the NSF-grants
DMS-0903856, DMS-0802927, and DMS-1047602. U.~Frauenfelder was partially supported by the Basic Research fund 2010�0007669 and O.~van Koert by the New Faculty Research Grant 0409-20100147 funded by the Korean government basic.

P.~Albers, J.~Fish, U.~Frauenfelder and O.~van Koert thank the IAS for its hospitality.

\section{History, known results and open questions}
Near the end of his lifelong quest to find periodic orbits, Poincar\'e introduced the concept of an annulus-like global surface of section (see \cite{Poincare}). In that same article, Poincar\'e observed that if a certain fixed point theorem (specifically Poincar\'e's last geometric theorem) holds true then the existence of such an annular surface of section implies the existence of periodic orbits. Shortly thereafter, Birkhoff proved Poincar\'e's last geometric theorem (see \cite{Birkhoff_Proof_of_Poincares_geometrict_heorem}) and then later generalized the notion of an annular surface of section to a surface of arbitrary genus and with an arbitrary number of boundary components (see \cite{Birkhoff_Dynamical_systems_with_two_degrees_of_freedom}).

The above results of Poincar\'e and Birkhoff were then employed by Conley in \cite{Conley_On_some_new_long_periodic_solutins_of_the_plane_restricted_three_body_problem} to prove the existence of certain long periodic orbits in the planar restricted three-body problem.  More precisely, Conley proved that there exists a sufficiently negative constant $E_0$, which is
independent of the mass ratio $\mu$, with the property that each energy surface $\{H=E<E_0\}$ admits an annulus-like surface of section. Under this assumption there are two bounded Hill's regions, and regularizing the associated singularities gives bounded energy surfaces. Each component of these sufficiently negative energy levels is diffeomorphic to $\R P^3$, and it is heuristically clear that they are well modeled by a small perturbation of the regularized Kepler problem. It is then possible to construct surfaces of section for the regularized Kepler problem, which persist under small perturbations.  Also note that an alternative approach using only canonical transformations can be found in \cite{Kummer_On_the_stability_of_Hills_solutions_of_the_plane_restricted_three_body_problem}.

For sufficiently small mass ratio McGehee \cite{McGehee_PhD} constructs disk-like surfaces of section around the heavy primary for energies up to the energy of the first Lagrange point. He also uses the Levi-Civita regularization and works in the double covering, $S^3$,  as we do. The purpose of this article is to prove the analogue of McGehee's theorem around the small primary. The surfaces of section in the articles by Conley and McGehee are perturbations of surfaces of section of the Kepler problem, which is completely integrable. We apply a result obtained by an entirely different method due to Hofer-Wysocki-Zehnder \cite{HWZ_the_dynamics_on_three_dimensional_strictly_conve_eneergy_surfaces}, which is based on holomorphic curve techniques. For this it suffices to prove convexity of the Levi-Civita embedding of the energy hypersurface into $\C^2$.

We would like to raise the following question:
\begin{question}  Does there exist a global (disk-like)  surface of section for each mass ratio $\mu$ and energy below the critical value $H(L_1)$
in both bounded energy components of the regularized problem?
\end{question}
The above question could be answered in the affirmative provided one could show that in the appropriate energy range
the two bounded components of the regularized energy surfaces are dynamically convex, i.e.~all Conley-Zehnder indices are greater or equal 3, see \cite{HWZ_the_dynamics_on_three_dimensional_strictly_conve_eneergy_surfaces}. There is an interesting recent paper
by U. Hryniewicz and P. Salomao, which gives a sufficient and necessary condition for the existence of global disk-like surfaces of section
which is relevant to our question, see \cite{Hryniewicz_Salomao_On_the_existence_of_disk_like_global_surfaces_of_section_for_Reeb_flows_on_the_tight_3_sphere}, and even goes beyond dynamical convexity.

It was shown in \cite{HWZ_the_dynamics_on_three_dimensional_strictly_conve_eneergy_surfaces} that strong convexity implies dynamical convexity. Also note that in Appendix \ref{sec:appendix} below we show that for energies near $H(L_1)$ the energy surface near the large primary fails to be convex. On the other hand we prove in \cite{Albers_Fish_Frauenfelder_Hofer_Koert_The_Conley_Zehnder_indices_of_the_rotating_Kepler_problem} that for the same energy levels, the energy surface near the primary at $0$ is \emph{dynamically} convex provided the mass ratio $\mu$ is sufficiently small;
this corresponds to a heavy primary located in $q=0$.
The method of the present paper can be used to check for a large class of pairs $(c,\mu)$ whether the energy hypersurface is indeed convex.

For energies just a bit higher than the critical value $H(L_1)$, the topology of (the bounded component of) the energy hypersurface changes from a disjoint union of two copies of $\R P^3$ to a connected sum $\R P^3\#\R P^3$. For topological reasons, global surfaces of section of disk or annulus type do not exist for  $\R P^3\#\R P^3$. However, the more general theory of finite energy foliations developed by Hofer-Wysocki-Zehnder in \cite{HWZ_Finite_energy_foliations_of_tight_three_spheres} still applies. We shall discuss this in an upcoming paper.

We expect the following global picture for energy levels just above $H(L_1)$. It is well-known that in the neck region of the connected sum there exists a hyperbolic periodic orbit with Conley-Zehnder index 2; this is the Lyapunov or halo orbit.
We expect there to exist a finite energy foliation where this Lyapunov orbit is one of at least three binding orbits. The existence of such a foliation would yield a structure theorem explaining in some detail the global behavior of the stable and unstable manifold of the Lyapunov orbit. Furthermore it would give a geometric explanation of the existence of a well-known homoclinic orbit asymptotic to this Lyapunov orbit, see Conley \cite{Conley_Twist_mappins_linking_analyticity} and McGehee \cite{McGehee_PhD}.
The reader should consult  \cite[Theorem 1.9]{HWZ_Finite_energy_foliations_of_tight_three_spheres}, where the theory of finite energy surfaces is developed for contact-type
flows on $S^3$. It is possible to make this technology work also for the connected sum of ${\mathbb RP}^3$'s.

Finally we would like to point out that convex energy surfaces have interesting symplectic  and dynamical properties. For example, the smallest occurring action of a periodic orbit
on a strongly convex energy surface is a symplectic capacity, which in turn is a crucially important concept in symplectic geometry. Also, it is very likely that in the case that $\Sigma\subset {\mathbb R}^4$, the surface of section is bounded by a periodic orbit of smallest action. This has been an an open problem for quite some time.  It is not too difficult to find these orbits by minimization
of a dual action functional, a method which can be implemented numerically. The details of the latter remarks are explained in \cite{Hofer_Zehner_Book}.
The relevant numerical methods are described in \cite{Goeing_Jaeschke_PhD}.

\section{Convexity of the planar restricted 3-body problem}

We recall that $\Sigma_{\mu,c}$ is the compact component of $K_{\mu,c}^{-1}(0)$ corresponding to the component of the energy hypersurface $H^{-1}(c)$ around earth, for a fixed mass ratio $\mu\in [0,1]$.
\begin{proposition}\label{thm:convex}
For every $c>\frac32$ there exists $\mu_0=\mu_0(c)\in (0,1)$ such that for all $\mu_0<\mu<1$ the component $\Sigma_{\mu,c}$ of the energy hypersurface of $K_{\mu,c}^{-1}(0)$ bounds a strongly convex domain.
\end{proposition}

\begin{proof}
We first compute the Hessian of
\beq
K_{\mu,c}(v,u)=\frac12|u|^2+c|v|^2+2|v|^2\langle u,iv\rangle-\mu\Im(uv)-\frac{\mu|v|^2}{|2v^2-1|}-\frac{1-\mu}{2}\;.
\eeq
In order to do so we need some auxiliary computations and consider the Hessian of
\beq\nonumber
g(v)=\frac{1}{|2v^2-1|}\;.
\eeq
For this we first set
\beq\nonumber
f(v)=|2v^2-1|^2=(2v^2-1)(2\bar{v}^2-1)=4|v|^4-4\Re(v^2)+1\;.
\eeq
Then we see
\beq\nonumber
df(v)\hat{v}=16|v|^2 \langle v,\hat{v}\rangle-8\Re(v\hat{v})\;,
\eeq
and
\beq\nonumber
D^2f(v)[\hat{v},\hat{v}]=32\langle v,\hat{v}\rangle^2+16|v|^2|\hat{v}|^2-8\Re(\hat{v}^2)\;.
\eeq
Thus
\bea\nonumber
Dg(v)\hat{v}&=-\tfrac12f(v)^{-\frac32}df(v)\hat{v}\\
&=-\frac{1}{|2v^2-1|^3}\big(8|v|^2\langle v,\hat{v}\rangle -4\Re(v\hat{v})\big)\;,
\eea
and
\bea\nonumber
D^2g(v)[\hat{v},\hat{v}]&=\tfrac34f(v)^{-\frac52}\big(df(v)\hat{v}\big)^2-\frac12f(v)^{-\frac32}D^2f(v)[\hat{v},\hat{v}]\\
&=\frac{3}{4|2v^2-1|^5}\big(16|v|^2\langle v,\hat{v}\rangle -8\Re(v\hat{v})\big)^2\\
&-\frac{4}{|2v^2-1|^3}\big(4\langle v,\hat{v} \rangle^2+2|v|^2|\hat{v}|^2-\Re(\hat{v}^2)\big)\;.\\
\eea
Therefore, we compute
\bea\nonumber
D\bigg(\frac{|v|^2}{|2v^2-1|}\bigg)\hat{v}&=-\frac{|v|^2}{|2v^2-1|^3}\big(8|v|^2 \langle v,\hat{v}\rangle-4\Re(v\hat{v})\big)\\
&+ \frac{2\langle v,\hat{v}\rangle}{|2v^2-1|}\;,
\eea
and finally
\bea\nonumber
D^2\bigg(\frac{|v|^2}{|2v^2-1|}\bigg)[\hat{v},\hat{v}]&=\frac{3|v|^2}{4|2v^2-1|^5}\big(16|v|^2\langle v,\hat{v}\rangle -8\Re(v\hat{v})\big)^2\\
&-\frac{4|v|^2}{|2v^2-1|^3}\big(4\langle v,\hat{v}\rangle^2+2|v|^2|\hat{v}|^2-\Re(\hat{v}^2)\big)\\
&-\frac{2\langle v,\hat{v}\rangle}{|2v^2-1|^3}\big(8|v|^2\langle v,\hat{v}\rangle-4\Re(v\hat{v})\big)\\
&+ \frac{2|\hat{v}|^2}{|2v^2-1|}\\
&-\frac{2\langle v,\hat{v}\rangle}{|2v^2-1|^3}\big(8|v|^2\langle v,\hat{v}\rangle-4\Re(v\hat{v})\big)\;.
\eea

We simplify this to
\bea\nonumber
D^2\bigg(\frac{|v|^2}{|2v^2-1|}\bigg)[\hat{v},\hat{v}]&=\frac{48|v|^2}{|2v^2-1|^5}\big(2|v|^2\langle v,\hat{v}\rangle-\langle\bar{v},\hat{v}\rangle\big)^2\\
&-\frac{4|v|^2}{|2v^2-1|^3}\big(4\langle v,\hat{v}\rangle^2+2|v|^2|\hat{v}|^2-\Re(\hat{v}^2)\big)\\
&-\frac{16\langle v,\hat{v}\rangle }{|2v^2-1|^3}\big(2|v|^2\langle v,\hat{v}\rangle -\langle\bar{v},\hat{v}\rangle\big)\\
&+ \frac{2|\hat{v}|^2}{|2v^2-1|}\;.
\eea

From this we conclude that the Hessian of
\beq\nonumber
K_{\mu,c}(v,u)=\frac12|u|^2+c|v|^2+2|v|^2\langle u,iv\rangle -\mu\Im(uv)-\frac{\mu|v|^2}{|2v^2-1|}-\frac{1-\mu}{2}
\eeq
is
\bea\nonumber
D^2K_{\mu,c}(u,v)[(\hat{u},\hat{v}),(\hat{u},\hat{v})]&=|\hat{u}|^2+2c|\hat{v}|^2+4\langle u,iv\rangle|\hat{v}|^2+8\langle v,\hat{v}\rangle\langle u,i\hat{v}\rangle\\
&+8\langle v,\hat{v}\rangle \langle \hat{u},iv\rangle+4|v|^2\langle \hat{u},i\hat{v}\rangle\\
&-2\mu\Im (\hat{u}\hat{v})\\
&-\frac{48\mu|v|^2}{|2v^2-1|^5}\big(2|v|^2\langle v,\hat{v}\rangle -\langle \bar{v},\hat{v}\rangle \big)^2\\
&+\frac{4\mu|v|^2}{|2v^2-1|^3}\big(4\langle v,\hat{v}\rangle ^2+2|v|^2|\hat{v}|^2-\Re(\hat{v}^2)\big)\\
&+\frac{16\mu\langle v,\hat{v}\rangle }{|2v^2-1|^3}\Big(2|v|^2\langle v,\hat{v}\rangle -\langle \bar{v},\hat{v}\rangle \Big)\\
&-\frac{2\mu|\hat{v}|^2}{|2v^2-1|}\;.
\eea

From now on we fix $c>\frac32$. We observe that in the limit $\mu\to1$ the energy hypersurface $\Sigma_{\mu,c}$ collapses onto the origin. To see this, first note that as $\mu\to 1$, the distance of $q$ to the Lagrange point goes to $0$.
As was also observed in \cite{Albers_Frauenfelder_Koert_Paternain_Liouville_field_for_PCR2BP}, this distance provides an upper bound for the size of Hill's region in the projection to the $(q_1,q_2)$-plane.
Since $q=2v^2$, we see that $|v|\to 0$ as $\mu\to 1$.
Now observe that the level set $K_{\mu,c}=0$ may be written as
$$
0=\frac12|u|^2+|u|\left(|v|^2\langle \frac{u}{|u|},iv \rangle -\mu\Im(\frac{u}{|u|} v ) \right )
+|v|^2 \left(c-\frac{\mu}{|1-2v^2|} \right)
-\frac{1-\mu}{2}\;.
$$
Regard this as a quadratic equation for $|u|$ which we can solve explicitly,
\bean
|u|=&-\left(|v|^2\langle \frac{u}{|u|},iv \rangle -\mu\Im(\frac{u}{|u|} v ) \right )\\[1ex]
&+
\sqrt{
\left(|v|^2\langle \frac{u}{|u|},iv \rangle -\mu\Im(\frac{u}{|u|} v ) \right)^2-2\left ( |v|^2(c-\frac{\mu}{|1-2v^2|})
-\frac{1-\mu}{2}\right )
}
\eea
Since $|v|\to 0$ and $1-\mu\to 0$ as $\mu\to 1$, we see that $|u|\to 0$ as $\mu\to 1$. In other words, $|u|,|v|\to 0$ as $\mu\to 1$ as claimed.  Consequently, for each $0<\epsilon<\frac14$ there exists $\mu_1=\mu_1(\epsilon)<1$ such that for $\mu_1\leq\mu\leq1$ we have
\beq\nonumber
|u|^2,\;|v|^2<\epsilon\quad\text{for all } (u,v)\in\Sigma_{\mu,c}\;.
\eeq
Thus, there exists a constant $C>3$ independent of $\epsilon$, $\mu$, and $c$ (provided that $c>\frac32$), such that
\beq\nonumber
D^2K_{\mu,c}(u,v)[(\hat{u},\hat{v}),(\hat{u},\hat{v})]\geq|\hat{u}|^2+2c|\hat{v}|^2-C\epsilon\big(|\hat{u}|^2+|\hat{v}|^2\big)-2\mu|\Im (\hat{u}\hat{v})|-\frac{2\mu|\hat{v}|^2}{|2v^2-1|}\;.
\eeq
Using again that $|v|^2<\epsilon<\tfrac14$ and the inequality $\frac{1}{1-x}\leq 1+2x$ for $0\leq x\leq\tfrac12$ we can estimate
\beq\nonumber
\frac{1}{|2v^2-1|}\leq1+4\epsilon\;,
\eeq
and thus
\bea\nonumber
D^2K_{\mu,c}(u,v)[(\hat{u},\hat{v}),(\hat{u},\hat{v})]&\geq|\hat{u}|^2+2c|\hat{v}|^2-C\epsilon\big(|\hat{u}|^2+|\hat{v}|^2\big)-2\mu|\Im (\hat{u}\hat{v})|-2\mu|\hat{v}|^2(1+4\epsilon)\;.
\eea
Estimating further, we note that
\beq\nonumber
2|\Im(\hat{u}\hat{v})|\leq\Delta| \hat{u} |^2+\frac1\Delta| \hat{v} |^2
\eeq
for each $\Delta>0$, and thus
\bea\nonumber
D^2K_{\mu,c}(u,v)[(\hat{u},\hat{v}),(\hat{u},\hat{v})]&\geq|\hat{u}|^2+2c|\hat{v}|^2-C\epsilon\big(|\hat{u}|^2+|\hat{v}|^2\big)-\mu(\Delta |\hat{u}|^2+\frac1\Delta|\hat{v}|^2)\\
&\phantom{\geq}-2\mu|\hat{v}|^2(1+4\epsilon)\\
&\geq|\hat{u}|^2(1-C\epsilon-\mu\Delta)+|\hat{v}|^2(2c-C\epsilon-\frac{\mu}{\Delta}-2\mu(1+4\epsilon))\;.
\eea
Next choose
\beq\nonumber
\epsilon<\min\left\{\frac{2c-3}{3C-8},\frac{1}{2C}\right\}\;,
\eeq
and fix $\Delta=1-C\epsilon$. Then by using $\mu<1$ we find
\bea\nonumber
D^2K_{\mu,c}(u,v)[(\hat{u},\hat{v}),(\hat{u},\hat{v})]&\geq|\hat{u}|^2(1-C\epsilon-\mu(1-C\epsilon))+|\hat{v}|^2(2c-C\epsilon-\frac{1}{1-C\epsilon}-2(1+4\epsilon))\\
&\geq|\hat{u}|^2(1-\mu)(1-C\epsilon)+|\hat{v}|^2(2c-2-(C+8)\epsilon-\frac{1}{1-C\epsilon})\;.
\eea
Using $\epsilon<\frac{1}{2C}$ we estimate as above
\beq\nonumber
\frac{1}{1-C\epsilon}\leq(1+2C\epsilon)\;,
\eeq
and thus
\bea\nonumber
D^2K_{\mu,c}(u,v)[(\hat{u},\hat{v}),(\hat{u},\hat{v})]&\geq|\hat{u}|^2(1-\mu)(1-C\epsilon)+|\hat{v}|^2(2c-2-(C+8)\epsilon-\frac{1}{1-C\epsilon})\\
&\geq|\hat{u}|^2(1-\mu)(1-C\epsilon)+|\hat{v}|^2(2c-2-(C+8)\epsilon-(1+2C\epsilon))\\
&=|\hat{u}|^2(1-\mu)(1-C\epsilon)+|\hat{v}|^2(2c-3-(3C+8)\epsilon)\;.
\eea
Finally since $\epsilon<\frac{2c-3}{3C-8}$ we obtain
\bea\nonumber
D^2K_{\mu,c}(u,v)[(\hat{u},\hat{v}),(\hat{u},\hat{v})]&\geq|\hat{u}|^2\underbrace{(1-\mu)(1-C\epsilon)}_{>0}+|\hat{v}|^2\underbrace{(2c-3-(3C+8)\epsilon)}_{>0}\;.
\eea
Thus, for a suitable $\delta=\delta(\mu,c)>0$ we have the desired estimate
\bea\nonumber
D^2K_{\mu,c}(u,v)[(\hat{u},\hat{v}),(\hat{u},\hat{v})]\geq \delta\cdot |(\hat{u},\hat{v})|^2\;.
\eea
In particular, the set $\Sigma_{\mu,c}$ bounds a strongly convex domain, which completes the proof of the proposition.
\end{proof}

\begin{proof}[Proof of Theorem \ref{thm:main}]
This follows from \cite[Theorem 1.3]{HWZ_the_dynamics_on_three_dimensional_strictly_conve_eneergy_surfaces} and Proposition \ref{thm:convex}.
\end{proof}

\appendix

\section{Convexity fails for high energy}\label{sec:appendix}

\begin{figure}[h]
\includegraphics[width=0.5\textwidth,clip]{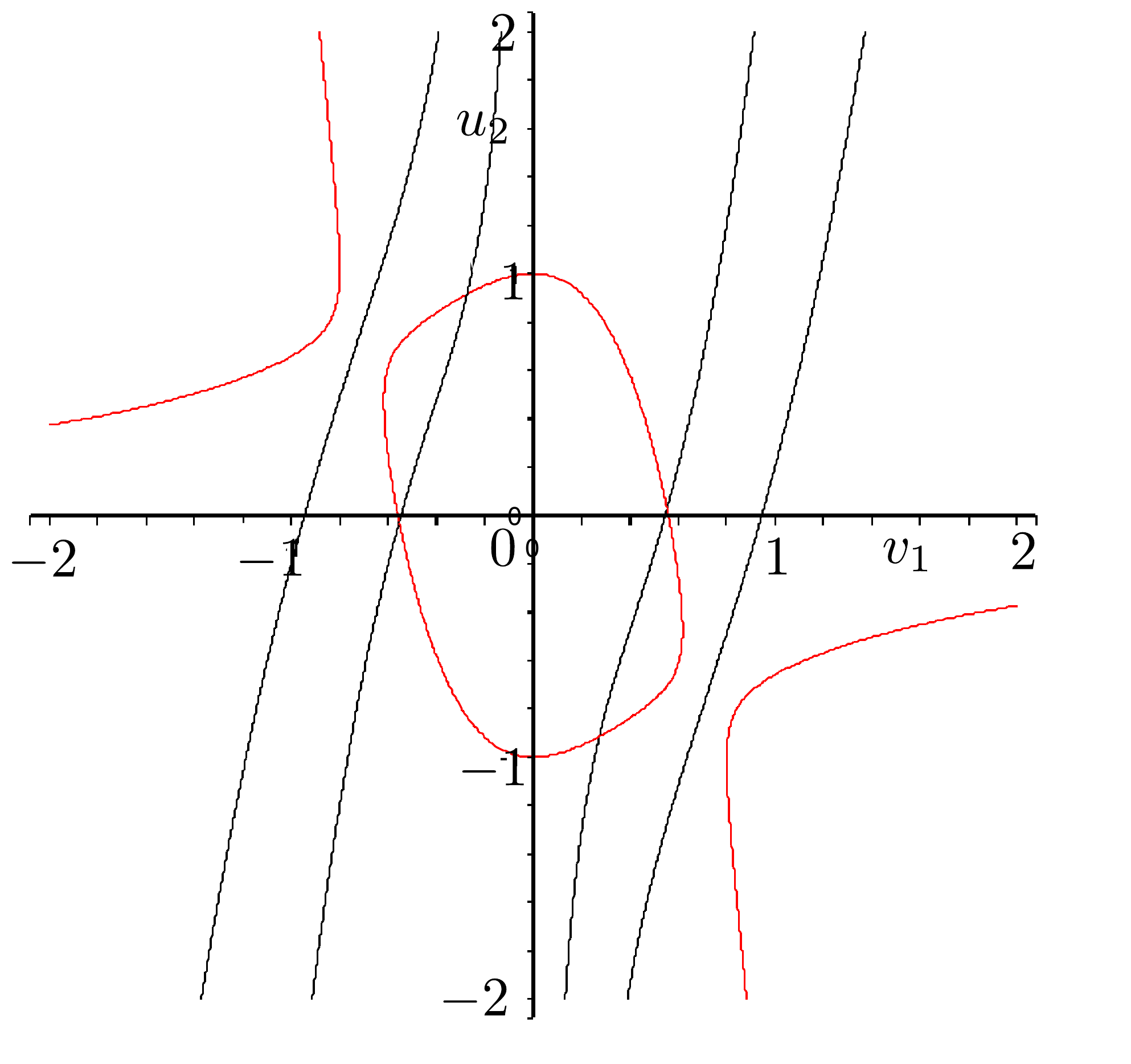}
\caption{The Levi-Civita embedding is not always convex: plot for $c=1.601$}\label{fig:not_convex}
\end{figure}

Here we shall briefly demonstrate that the Levi-Civita embedding fails to be convex for energies close to the Lagrange energy.
In Figure \ref{fig:not_convex} we shall plot the level set $\{ K=0 \}$ intersected with the hyperplanes $v_2=0$ and $u_1=0$.
In the same picture we also plot where $\det  D^2K=0$.
Since the curves intersect, we see that convexity of the level set $\{K=0\}$ fails at those points.

The computations in this section were done with MAPLE.
As often is the case, expressions obtained by computer algebra programs tend to be rather unwieldy, but in case $\mu=0$, the expression for $\det D^2 K$ is still manageable,
\begin{align*}
\det D^2 K&=
2304\,v_1^{8}+9216\,v_1^{6}v_2^{2}-3072\,u_2 v_1^{5}+13824\,v_1^{4}v_2^{4}-1280\,cv_1^{
4}+3072\, u_1 v_1^{4} v_2\\
&\qquad-6144\, u_2 v_1^{3} v_2^{2}+6144\,u_1 v_1^{2} v_2^{3}-2560\,cv_1^{2}v_2^{2}-256\,u_1^{2} v_1^{2}+9216 \,v_1^{2}v_2^{6}\\
&\qquad+768\,u_2^{2} v_1^{2}-3072 \,u_2 v_1 v_2^{4}+512\,c u_2 v_1 -2048\,u_1 u_2 v_1  v_2 +64\,c^{2}\\
&\qquad+2304\,v_2^{8} +768\,u_1^{2} v_2^{2}-512\,c u_1 v_2 +3072\,u_1 v_2^{5}-1280\,c v_2^{4}-256\,u_2^{2} v_2^{2}.
\end{align*}
%\clearpage

\bibliographystyle{amsalpha}
\bibliography{../../../../../Bibtex/bibtex_paper_list}

%\begin{thebibliography}{999}
%
%%%%%% \bibitem{AFKP}
%%%%% P.\,Albers, U.\,Frauenfelder, O.\,van Koert, G.\,Paternain, {\em
%%%%%%%%%%% The contact geometry of the restricted 3-body problem},
%%% arXiv:1012.2140.
%\bibitem{CFK} K.\,Cieliebak, U.\,Frauenfelder, O.\,van Koert,
%{\em The Cartan geometry of the rotating Kepler problem}, preprint.
%\bibitem{D} J.\,Duistermaat, {\em On the Morse index in variational
%calculus}, Advances in Math.\,\textbf{21}:2 (1976), 173--195.
%\bibitem{M} J.~Moser, {\em Regularization of Kepler's problem and the
%averaging method on a manifold}, Comm.\,Pure Appl.\,Math.\textbf{23}
%(1970) 609-636.
%\bibitem{W} J.\,Weber, {\em Perturbed closed geodesics are periodic
%orbits:\,index and transversality}, Math.\,Z.\,\textbf{241}:1
%(2002), 45--82.
%\end{thebibliography}

\end{document}